\documentclass[12pt]{amsart}
\usepackage{amsfonts}
\usepackage{amssymb}
\theoremstyle{definition}

\theoremstyle{remark}

\newcommand{\ds}{\displaystyle}

\begin{document}

\centerline{\bf Comptes rendus de l'Acad\'emie bulgare des Sciences}

\centerline{\it Tome 35, No 3, 1982}

\vspace{0.6in}

\begin{flushright}
{\it MATH\'EMATIQUES
\\ G\'eometrie diff\'erentielle}
\end{flushright}

\vspace{0.2in}

\centerline{\large\bf ON THE AXIOM OF SPHERES IN K\"AHLER GEOMETRY }

\vspace{0.3in}
\centerline{\bf O. T. Kassabov}

\vspace{0.3in}
\centerline{\it (Submited by Academician B. Petkanchin on November 25, 1981)}

\vspace{0.2in}

{\bf 1. Introduction}. Let $M$ be a $2m$-dimensional K\"ahler manifold with a Riemannian 
metric $g$ and a complex structure $J$. Let $R,\,S $ and $\tau(R)$ denote the curvature 
tensor, the Ricci tensor and the scalar curvature of $M$, respectively. The Bochner
curvature tensor $B$ is defined by
$$
	\begin{array}{c}
		B(X,Y,Z,U)=R(X,Y,Z,U) - \frac{1}{2(m+2)} \{ g(X,U)S(Y,Z)-g(X,Z)S(Y,U)  \\
		+g(Y,Z)S(X,U)-g(Y,U)S(X,Z) + g(X,JU)S(Y,JZ)-g(X,JZ)S(Y,JU)  \\
		+g(Y,JZ)S(X,JU)-g(Y,JU)S(X,JZ) - 2g(X,JY)S(Z,JU)-2g(Z,JU)S(X,JY) \}  \\
		+\frac{\tau(R)}{4(m+1)(m+2)} \{ g(X,U)g(Y,Z)-g(X,Z)g(Y,U) \\
		+g(X,JU)g(Y,JZ)-g(X,JZ)g(Y,JU) -2g(X,JY)g(Z,JU) \} \ .
	\end{array}
$$ 

Let $N$ be an $n$-dimensional submanifold of $M$. The second fundamental form $\alpha$
of the immersion is defined by 
$\alpha(X,Y)=\widetilde\nabla_XY - \nabla_XY$ for $X,\ Y \in {\mathfrak X}(N)$, where
$\widetilde\nabla$ (resp. $\nabla$) is the covariant differentiations on $M$ (resp. $N$).
The submanifold $N$ is said to be totally
umbilical, if $\alpha(X,Y)=g(X,Y)H$, $H$ being the mean curvature vector of $N$ in $M$,
i.e.  $H=(1/n){\rm trace}\, \alpha$. Let $\xi $ be a vector field normal to $N$. Then 
the Weingarten formula is
$$ 
	\widetilde\nabla_X\xi =- A_{\xi}X+D_X{\xi} \ ,
$$
where $-A_{\xi}X$ (respectively, $D_X\xi$)
denotes the tangential (resp. the normal) component of $\widetilde\nabla_X\xi$. The
vector field $\xi$ is said to be parallel, if $D_X\xi=0$ for each $X\in \mathfrak X(N)$.

The equation of Codazzi is given by
$$
	\{R(X,Y)Z\}^{\perp} = (\overline\nabla_X\alpha)(Y,Z)-(\overline\nabla_Y\alpha)(X,Z) \ ,  
$$
where $\{R(X,Y)Z\}^{\perp}$ denotes the normal component of $R(X,Y)Z$ and 
$$
	(\overline\nabla_X\alpha)(Y,Z)=D_X\alpha(Y,Z)-\alpha(\nabla_XY,Z)-\alpha(Y,\nabla_XZ) \ . 
$$

By an $n$-plane in $T_p(M)$ we mean an $n$-dimensional linear subspace of  $T_p(M)$. 
An $n$-plane  $\pi$  is said to be 
holomorphic (resp. antiholomorphic) if  \ $\pi =J \pi$ \ 
(resp. \ $\pi \perp J\pi)$.

A $2m$-dimensional K\"ahler manifold $M$ is said to satisfy the axiom of holomorphic 
(resp. antiholomorphic) $2n$-spheres (resp. $n$-spheres), where $n$ is a fixed integer,
$1\le n \le m$, if for each point $p \in M$ and for any $2n$-dimensional
holomorphic (resp. $n$-dimensional antiholomorphic) plane \ $\pi$ \ in \ $T_p(M)$ \ there exists a totally umbilical
submanifold \ $N$ of $M$  with a 
parallel mean curvature vector, such that \ $p\in N$ and \ $T_pN=\pi$.

As was proved in [2], if a $2m$-dimensional K\"ahler manifold $M$ satisfies
the axiom of holomorphic $2n$-spheres for some $n$, $1\le n<m$ or the axiom of
antiholomorphic  $n$-spheres for some $n$, $1< n\le m$, then $M$ is of constant
holomorphic sectional curvature.

We shall prove the following theorems:

\vspace{0.05in}
{\bf Theorem 1.} Let $M$ be a $2m$-dimensional K\"ahler manifold, $m>2$, and let
$n$ be a fixed integer, $2\le n <m$. If for each point $p\in M$ and for any 
holomorphic $2n$-plane $\pi$ in $T_p(M)$ there exists a totally umbilical submanifold 
$N$ of $M$, such that $p \in N$ and $T_p(N)=\pi$, then $M$ is of constant holomorphic 
sectional curvature.

\vspace{0.05in}
{\bf Theorem 2.} Let $M$ be a $2m$-dimensional K\"ahler manifold, $m>2$, and let
$n$ be a fixed integer, $2< n \le m$. If for each point $p\in M$ and for any 
antiholomorphic $n$-plane $\pi$ in $T_p(M)$ there exists a totally umbilical submanifold 
$N$ of $M$, such that $p \in N$ and $T_p(N)=\pi$, then $M$ is of constant holomorphic 
sectional curvature.

\vspace{0.2in}
{\bf 2. A Lemma and Proofs of the Theorems.} As is known, a $2m$-dimensional K\"ahler
manifold $M$ has vanishing Bochner curvature tensor, iff for each point $p$ of $M$
the sum $\ds \sum_{i=1}^m R(e_i,Je_i,Je_i,e_i)$ is independent of the orthonormal
basis $ \{ e_i,Je_i;\,i=1,...,m \}$ of $T_p(M)$ [5]. Hence it is not difficult to
prove the following

\vspace{0.05in}
{\bf Lemma.} A K\"ahler manifold $M$ of dimension $2m \ge 6$ has a vanishing Bochner
curvature tensor, iff for each point $p\in M$ and for all unit vectors 
$x,\,y,\,z \in T_p(M)$ which span an antiholomorphic 3-plane
$$
	R(x,Jx,y,z)=2R(x,y,Jx,z)
$$
holds good.

\vspace{0.05in}
Let $N$ be a totally umbilical submanifold of $M$. Then Codazzi's equation reduces to
$$
	\{ R(X,Y)Z\}^{\perp} = g(Y,Z)D_XH-g(X,Z)D_YH \ .  \leqno (2.1)
$$

Now we can proceed to prove Theorem 1. For a point $p \in M$ we take arbitrary
unit vectors $x,\,y,\,z \in T_p(M)$ which span an antiholomorphic 3-plane. Let $N$ be a
totally umbilical submanifold of $M$ such that $p \in N$, $x,\,y,\,Jx,\,Jy \in T_p(N)$
and $z \perp T_p(N)$. Then, from (2.1) we obtain
$$
	R(x,Jx,y,z)=0 \ ,       \leqno (2.2)
$$
$$
	R(x,y,Jx,z)=0
$$
and, according to the Lemma, $M$ has vanishing Bochner curvature tensor. Hence
$$
	\begin{array}{c}
		R(X,Y,Z,U) = \frac{1}{2(m+2)} \{ g(X,U)S(Y,Z)  \\
		-g(X,Z)S(Y,U)+g(Y,Z)S(X,U)-g(Y,U)S(X,Z)   \\
		+ g(X,JU)S(Y,JZ)-g(X,JZ)S(Y,JU)+g(Y,JZ)S(X,JU)  \\
		-g(Y,JU)S(X,JZ) - 2g(X,JY)S(Z,JU)-2g(Z,JU)S(X,JY) \}  \\
		-\frac{\tau(R)}{4(m+1)(m+2)} \{ g(X,U)g(Y,Z)-g(X,Z)g(Y,U) \\
		+g(X,JU)g(Y,JZ)-g(X,JZ)g(Y,JU) -2g(X,JY)g(Z,JU) \} \ .
	\end{array}      \leqno (2.3)
$$ 
From (2.2) and (2.3)
$$
	S(y,z)=0
$$
for all $y,\,z \in T_p(M)$ with $g(y,z)=g(y,Jz)=0$ and for each point $p \in M$. 
Consequently, $M$ is an Einsteinian manifold and because of $B=0$ \ $M$ is of
constant holomorphic sectional curvature. 

\vspace{0.05in}
The proof of Theorem 2 is similar.

\vspace{0.05in}
{\bf Remark 1.} If $M$ is a K\"ahler manifold of constant holomorphic sectional 
curvature, then every totally umbilical submanifold of $M$ has parallel mean 
curvature vector, see [3].

\vspace{0.05in}
{\bf Remark 2.} It is known that any $2m$-dimensional K\"ahler manifold of constant
holomor\-phic sectional curvature satisfies the axiom of holomorphic $2n$-spheres
and the axiom of antiholomorphic $n$-spheres for each $n$, $1\le n\le m$, see e.g. [3].

\vspace{0.05in}
{\bf Remark 3.} For a Riemannian manifold the condition analogous to that
in Theorem 1 or Theorem 2 is equivalent to the requirement that the manifold
has vanishing Weil curvature tensor, see [1].

\vspace {0.7cm}
\begin{flushright}
{\it Institute of Mathematics \\
Bulgarian Academy of Sciences \\
PO Box 373  \\
1090 Sofia, Bulgaria}
\end{flushright}

\end{document}